\documentclass[a4paper,12pt]{article}

\usepackage[english]{babel} 
\usepackage[ansinew]{inputenc}
\usepackage[LGR,T1]{fontenc}

\usepackage{verbatim}
\usepackage{color}
\usepackage{calc}
\usepackage{hyperref}
\usepackage[dvips]{graphicx}
\hypersetup{colorlinks=true, linkcolor=blue, filecolor=blue, urlcolor=blue}

\usepackage{amssymb}
\usepackage[mathscr]{eucal}
\usepackage{amscd,amsmath,amsfonts,amssymb,textcomp}
\usepackage[all]{xy}
\everymath{\displaystyle}

\newcommand{\MBC}{\mathbb{C}}

\newcommand{\MBN}{\mathbb{N}}

\newcommand{\MBQ}{\mathbb{Q}}

\newcommand{\MBZ}{\mathbb{Z}}

\newcommand{\MCA}{\mathcal{A}}
\newcommand{\MCB}{\mathcal{B}}
\newcommand{\MCC}{\mathcal{C}}

\newcommand{\MCE}{\mathcal{E}}

\newcommand{\MCG}{\mathcal{G}}
\newcommand{\MCH}{\mathcal{H}}
\newcommand{\MCI}{\mathcal{I}}
\newcommand{\MCJ}{\mathcal{J}}

\newcommand{\MCM}{\mathcal{M}}

\newcommand{\MCO}{\mathcal{O}}

\newcommand{\MCU}{\mathcal{U}}

\newcommand{\MCW}{\mathcal{W}}

\newcommand{\MFa}{\mathfrak{a}}

\newcommand{\MFf}{\mathfrak{f}}

\newcommand{\MFg}{\mathfrak{g}}

\newcommand{\MFh}{\mathfrak{h}}

\newcommand{\MFl}{\mathfrak{l}}

\newcommand{\MFm}{\mathfrak{m}}

\newcommand{\MFp}{\mathfrak{p}}

\newcommand{\MFq}{\mathfrak{q}}

\newcommand{\MFr}{\mathfrak{r}}

\newcommand{\MFu}{\mathfrak{u}}

\newcommand{\MFX}{\mathfrak{X}}

\newcommand{\MSC}{\mathscr{C}}

\newcommand{\GGa}{\alpha}
\newcommand{\GGb}{\beta}

\newcommand{\GVf}{\varphi}
\newcommand{\GGg}{\gamma}
\newcommand{\GGG}{\Gamma}
\newcommand{\GGi}{\iota}
\newcommand{\GGk}{\kappa}

\newcommand{\GGl}{\lambda}

\newcommand{\GGs}{\sigma}

\newcommand{\GGu}{\upsilon}

\newcommand{\GGz}{\zeta}

\newcommand{\dd}{\mathrm{d}}

\newcommand{\vid}{\varnothing}

\newcommand{\lA}{\left\{}

\newcommand{\la}{\left\langle}
\newcommand{\ra}{\right\rangle}

\newcommand{\alg}{\mathrm{alg}}
\newcommand{\Char}{\mathrm{char}}

\newcommand{\Cl}{\mathrm{Cl}}

\newcommand{\Gal}{\mathrm{Gal}}

\newcommand{\IM}{\mathrm{Im}}
\newcommand{\Ker}{\mathrm{Ker}}

\newcommand{\nr}{\mathrm{nr}}

\newcommand{\PR}{\mathrm{pr}}

\usepackage{makeidx}
\title{INVARIANTS AND COINVARIANTS OF SEMILOCAL UNITS MODULO ELLIPTIC UNITS.}
\author{\textsc{St\'ephane VIGUI\'E}
\footnote{St\'ephane Vigui\'e, Laboratoire de math\'ematiques de Besan\c con, UMR CNRS 6623, Universit\'e de Franche-Comt\'e, 16 route de Gray, 25030 Besançon cedex, France.
e-mail: \texttt{stephane.viguie@univ-fcomte.fr}}
}
\makeindex

\newtheorem{dfe}{dfe}[section]
\newtheorem{lem}{lem}[section]
\newtheorem{pro}{pro}[section]
\newtheorem{rem}{rem}[section]
\newtheorem{teh}{teh}[section]

\newtheorem{df}[dfe]{Definition}
\newtheorem{lm}[lem]{Lemma}
\newtheorem{pr}[pro]{Proposition}

\newtheorem{rmq}[rem]{Remark}
\newtheorem{theo}[teh]{Theorem}

\numberwithin{equation}{section}

\begin{document}
\maketitle

\begin{abstract}
Let $p$ be a prime number, and let $k$ be an imaginary quadratic number field in which $p$ decomposes into two primes $\MFp$ and $\bar{\MFp}$.
Let $k_\infty$ be the unique $\MBZ_p$-extension of $k$ which is unramified outside of $\MFp$, and let $K_\infty$ be a finite extension of $k_\infty$, abelian over $k$.
Let $\MCU_\infty/\MCC_\infty$ be the projective limit of principal semi-local units modulo elliptic units.
We prove that the various modules of invariants and coinvarants of $\MCU_\infty/\MCC_\infty$ are finite.
Our approach uses distributions and the $p$-adic $\mathrm{L}$-function, as defined in \cite{deshalit87}.
\end{abstract}

\noindent{\small\textbf{Mathematics Subject Classification (2010):} 11G16, 11R23.}
\\

\noindent{\small\textbf{Key words:} Elliptic units, Iwasawa theory.}

\section{Introduction.}

Let $p$ be a prime number, and let $k$ be an imaginary quadratic number field in which $p$ decomposes into two distinct primes $\MFp$ and $\bar{\MFp}$.
Let $k_\infty$ be the unique $\MBZ_p$-extension of $k$ which is unramified outside of $\MFp$, and let $K_\infty$ be a finite extension of $k_\infty$, abelian over $k$.
Let $G_\infty$ be the Galois group of $K_\infty/k$.
We choose a decomposition of $G_\infty$ as a direct product of a finite group $G$ (the torsion subgroup of $G_\infty$) and a topological group $\GGG$ isomorphic to $\MBZ_p$, $G_\infty = G\times \GGG$.
For all $n\in\MBN$, let $K_n$ be the field fixed by $\GGG_n := \GGG^{p^n}$, and let $G_n:=\Gal\left(K_n/k\right)$.
Remark that there may be different choices for $\GGG$, but when $p^n$ is larger than the order of the $p$-part of $G$, the group $\GGG_n$ does not depend on the choice of $\GGG$.

Let $F/k$ be a finite abelian extension of $k$.
We denote by $\MCO_F$ the ring of integers of $F$.
Then we write $\MCO_F^\times$ for the group of global units of $F$, and $C_F$ for the group of elliptic units of $F$ (see section \ref{sectionellipticunits}).
We set $\MCC_F:=\MBZ_p\otimes_\MBZ C_F$.
For all prime ideal $\MFq$ of $\MCO_F$ above $\MFp$, we write $F_\MFq$, $\MCO_{F_\MFq}$, and $\MCO_{F_\MFq}^\times$ respectively for the completion of $F$ at $\MFq$, the ring of integers  of $F_\MFq$, and the group of units of $\MCO_{F_\MFq}$.
Then we write $\MCU_F$ for the pro-$p$-completion of $\prod_{\MFq\vert\MFp} \MCO_{F_\MFq}^\times$.
The injection $\MCO_F^\times \hookrightarrow \prod_{\MFq\vert\MFp} \MCO_{F_\MFq}^\times$ induces a canonical map $\MBZ_p\otimes_\MBZ\MCO_F^\times \rightarrow \MCU_F$.
The Leopoldt conjecture, which is known to be true for abelian extensions of $k$, states that this map is injective.
For all $n\in\MBN$, we write $\MCC_n$ and $\MCU_n$ for $\MCC_{K_n}$ and $\MCU_{K_n}$ respectively.
We define $\MCC_\infty := \varprojlim\, \MCC_n$ and $\MCU_\infty := \varprojlim \, \MCU_n$ by taking projective limit under the norm maps.
The injections $\MCC_n\hookrightarrow\MCU_n$ are norm compatibles and taking the limit we obtain an injection $\MCC_\infty \hookrightarrow \MCU_\infty$.

For any profinite group $\MCG$, and any commutative ring $R$, we define the Iwasawa algebra
\[R\left[\left[\MCG\right]\right] := \varprojlim \, R\left[\MCH\right],\]
where the projective limit is over all finite quotient $\MCH$ of $\MCG$.
Then $\MCC_\infty$ and $\MCU_\infty$ are naturally $\MBZ_p\left[\left[G_\infty\right]\right]$-modules.
It is well known that they are finitely generated over $\MBZ_p[[\GGG]]$. 
Moreover one can show that $\MCU_\infty / \MCC_\infty$ is torsion over $\MBZ_p[[\GGG]]$ (see \cite[Proposition 3.1]{viguie11b}).
Let us fix a topological generator $\GGg$ of $\GGG$, and set $T:=\GGg-1$.
We denote by $\MBC_p$ a completion of an algebraic closure of $\MBQ_p$.
For any complete subfield $L$ of $\MBC_p$, finitely ramified over $\MBQ_p$, we denote by $\MCO_L$ the complete discrete valuation ring of integers of $L$.
Then the ring $\MCO_L[[\GGG]]$ is isomorphic to $\MCO_L[[T]]$.
It is well known that $\MCO_L[[T]]$ is a noetherian, regular, local domain.
We also recall that $\MCO_L[[T]]$ is a unique factorization domain.
If $\MFu_L$ is a uniformizer of $\MCO_L$, then the maximal ideal $\MFm_L$ of $\MCO_L$ is generated by $\MFu_L$ and $T$, and $\MCO_L[[T]]$ is a complete topological ring with respect to its $\MFm_L$-adic topology.
A morphism $f:M\rightarrow N$ between two finitely generated $\MCO_L[[T]]$-module is called a pseudo-isomorphism if its kernel and its cokernel are finitely generated and torsion over $\MCO_L$.
If a finitely generated $\MCO_L[[T]]$-module $M$ is given, then one may find elements $P_1$, ..., $P_r$ in $\MCO_L[T]$, irreducible in $\MCO_L[[T]]$, and nonnegative integers $n_0$, ..., $n_r$, such that there is a pseudo-isomorphism
\[M\longrightarrow \MCO_L[[T]]^{n_0} \oplus \bigoplus_{i=1}^r \MCO_L[[T]]/\left(P_i^{n_i}\right).\]
Moreover, the integer $n_0$ and the ideals $\left(P_1^{n_1}\right)$, ..., $\left(P_r^{n_r}\right)$, are uniquely determined by $M$.
If $n_0=0$, then the ideal generated by $P_1^{n_1}\cdots P_r^{n_r}$ is called the characteristic ideal of $M$, and is denoted by $\Char_{\MCO_L[[T]]}(M)$.

Let $\chi$ be an irreducible $\MBC_p$-character of $G$.
Let $L(\chi)\subset\MBC_p$ be the abelian extension of $L$ generated by the values of $\chi$.
The group $G$ acts naturally on $L(\chi)$ if we set, for all $g\in G$ and all $x\in L(\chi)$, $g.x:=\chi(g)x$.
For any $\MCO_L[G]$-module $Y$, we define the $\chi$-quotient $Y_\chi$ by $Y_\chi:=\MCO_{L(\chi)}\otimes_{\MCO_L[G]}Y$.
If $Y$ is an $\MCO_L\left[\left[G_\infty\right]\right]$-module, then $Y_\chi$ is an $\MCO_{L(\chi)}[[T]]$-module in a natural way.
Moreover if $L$ contains a $\left[K_0:k\right]$-th primitive root of unity, then there is $(a,b)\in\MBN^2$ such that
\begin{equation}
\label{chipart}
\MFu_L^a\Char_{\MCO_L[[T]]}(M) \, = \, \MFu_L^b\prod_{\chi} \Char_{\MCO_L[[T]]} \left(M_\chi\right),
\end{equation}
where the product is over all irreducible $\MBC_p$-character on $G$.

For any profinite group $\MCG$, any normal subgroup $\MCH$ of $\MCG$ and any $\MCO_L\left[\left[\MCG\right]\right]$-module $M$, we denote by $M^\MCH$ the module of $\MCH$-invariant of $M$, that is to say the maximal submodule of $M$ which is invariant under the action of $\MCH$.
We denote by $M_\MCH$ the module of $\MCH$-coinvariant of $M$, which is the quotient of $M$ by the closed submodule topologically generated by the elements $(h-1)m$ with $h\in\MCH$ and $m\in M$.

In this article, we prove that for all $n\in\MBN$, the module of $\GGG_n$-invariants and the module of $\GGG_n$-coinvariants of $\MCU_\infty / \MCC_\infty$ are finite (see Theorem \ref{MCUMCCinvcoinvfinite}).
It generalizes a part of a result of Coates-Wiles \cite[Theorem 1]{coates-wiles78}, where this result is shown at the $\chi^i$-parts, for $i\not\equiv0$ modulo $p-1$, and for $\chi$ the character giving the action of $G$ on the $\MFp$-torsion points of a certain elliptic curve.
But the result of \cite{coates-wiles78} is stated for non-exceptional primes $p$ (in particular $p\notin\{2,3\}$), and under the assumption that $\MCO_k$ is principal.
Here we prove the general case.

Moreover we would like to mention an application of Theorem \ref{MCUMCCinvcoinvfinite} to the main conjecture of Iwasawa theory.
For all $n\in\MBN$, we set $\MCE_n := \MBZ_p\otimes_\MBZ \MCO_{K_n}^\times$ and we denote by $A_n$ the $p$-part of the class-group $\Cl\left(\MCO_{K_n}\right)$ of $\MCO_{K_n}$.
We define $\MCE_\infty := \varprojlim\,\MCE_n$ and $A_\infty := \varprojlim\,A_n$, projective limits under the norm maps.
A formulation of the (one variable) main conjecture says that $\Char_{\MBZ_p(\chi)[[T]]}\left(\MCE_\infty/\MCC_\infty\right)_\chi = \Char_{\MBZ_p(\chi)[[T]]}\left(A_{\infty,\chi}\right)$,
where $\MBZ_p(\chi)$ is the ring of integers of $\MBQ_p(\chi)$.
It has been proved in many cases by the use of Euler systems.
We refer the reader to the pioneering work of Rubin in \cite[Theorem 4.1]{rubin91} and \cite[Theorem 2]{rubin94}, adapted to the cyclotomic case by Greither in \cite[Theorem 3.2]{greither92}.
The method is now classical, applied by many authors, see \cite[Theorem 3.1]{bley06}, \cite{oukhaba10} and \cite{viguie11b}. 
However for $p=2$ (and for some cases for $p=3$) the results are less strong and we just obtain a divisibility relation
\begin{equation}
\Char_{\MBZ_p(\chi)[[T]]}\left(A_{\infty,\chi}\right) \quad\text{divides}\quad p^a\Char_{\MBZ_p(\chi)[[T]]}\left(\MCE_\infty/\MCC_\infty\right)_\chi,
\label{divisibilityeh}
\end{equation}
for some $a\in\MBN$ (see \cite{oukhaba10} and \cite{viguie11b}).
Following the ideas of Belliard in \cite{belliard09}, in a forthcoming paper we will deduce from Theorem \ref{MCUMCCinvcoinvfinite} that for $p\in\{2,3\}$ the $\MBZ_p[[\GGG]]$-modules $\MCE_\infty / \MCC_\infty$ and $A_\infty$ have the same Iwasawa's $\mu$ and $\GGl$ invariants.
This result, together with (\ref{divisibilityeh}), implies that there is $(a,b)\in\MBN^2$ such that the following raw form of the main conjecture holds,
\[\MFu_\chi^a \Char_{\MBZ_p(\chi)[[T]]} \left(A_{\infty,\chi}\right) = \MFu_\chi^b\Char_{\MBZ_p(\chi)[[T]]} \left(\MCE_\infty/\MCC_\infty\right)_\chi,\]
where $\MFu_\chi$ is a uniformizer of $\MBZ_p(\chi)$.

\section{Distributions.}

In this section, let $A$ be a commutative ring and let $\MCG$ be a profinite group.
We denote by $\MFX(\MCG)$ the set of compact-open subsets of $\MCG$.
Remark that for any $X\in\MFX(\MCG)$, one can find a finite subset $F$ of $X$, and an open normal subgroup $\MCH$ of $\MCG$, such that $X=\mathop{\cup}_{x\in F}x\MCH$.

\begin{df}
An $A$-distribution on $\MCG$ is an application $\mu:\MFX(\MCG)\rightarrow A$, such that for all $(X_1,X_2)\in\MFX(\MCG)^2$, if $X_1\cap X_2=\vid$, then
\[ \mu\left(X_1\cup X_2\right) = \mu(X_1) + \mu(X_2).\]
We denote by $\MCM\left(\MCG,A\right)$ the $A$-module of $A$-distributions on $\MCG$.
Moreover for $X\in\MFX(\MCG)$ and $\mu\in\MCM\left(\MCG,A\right)$, we denote by $\mu_{\vert X}$ the restriction of $\mu$ to $X$, 
\[\mu_{\vert X} : \MFX(X) \xymatrix{\ar[r] &} A,\quad Y \xymatrix{\ar@{|->}[r] &} \mu(Y).\]
\end{df}

Let $\pi:\MCG\rightarrow\MCG'$ be a continuous open morphism between two profinite groups, such that $\Ker(\pi)$ is finite.
To any distribution $\mu\in\MCM\left(\MCG,A\right)$ we attach the unique $A$-distribution $\pi_\ast\mu$ on $\MCG'$, such that for all $X\in\MFX\left(\MCG'\right)$, 
\[\pi_\ast\mu(X) = \mu\left(\pi^{-1}(X)\right).\]
For any $\GGs\in\MCG$, let us also denote by $\GGs_\ast\mu$ the unique $A$-distribution on $\MCG$, such that for all $X\in\MFX(\MCG)$, 
\[\GGs_\ast\mu(X) = \mu\left(\GGs^{-1}X\right).\]
To any distribution $\mu'\in\MCM\left(\MCG',A\right)$ we attach the unique $A$-distribution $\pi^\sharp\mu'$ on $\MCG$, such that for all $g\in\MCG$, and all open subgroup $\MCH$ of $\MCG$,
\begin{equation}
\label{defpibullet} 
\pi^\sharp\mu' \left(g\MCH\right) \; = \; \#\left(\MCH\cap\Ker(\pi)\right) \, \mu' \left( \pi\left(g\MCH\right) \right). 
\end{equation}
Then we have
\begin{equation}
\label{pibulletpiast}
\pi^\sharp\pi_\ast\mu \; = \; \sum_{\GGs\in\Ker(\pi)} \GGs_\ast\mu \quad\text{and}\quad \pi_\ast\pi^\sharp\mu' \; = \; \#\left(\Ker(\pi)\right) \, \mu'_{\vert \IM(\pi)}.
\end{equation}
For $(\GGa_1,\GGa_2)\in\MCM\left(\MCG,A\right)^2$, there is a unique $A$-distribution $\GGb$ on $\MCG\times\MCG$ such that for all $(X_1,X_2)\in\MFX(\MCG)^2$, $\GGb\left(X_1\times X_2\right) = \GGa_1(X_1)\GGa_2(X_2)$.
Then the convolution product $\GGa_1\GGa_2$ of $\GGa_1$ and $\GGa_2$ is defined by 
$\GGa_1\GGa_2 := m_\ast\GGb$, 
where $m:\MCG\times\MCG\rightarrow\MCG$, $(\GGs_1,\GGs_2)\mapsto\GGs_1\GGs_2$.
Once equipped with the convolution product, $\MCM\left(\MCG,A\right)$ is an $A$-algebra.
For any $A$-distribution $\mu$ on $\MCG$, let us denote by $\underline{\mu}$ the unique element in $A\left[\left[\MCG\right]\right]$ such that for all open normal subgroup $\MCH$ of $\MCG$, the image $\underline{\mu}_\MCH$ of $\underline{\mu}$ in $A\left[\MCG/\MCH\right]$ is given by
\[\underline{\mu}_\MCH \, = \, \sum_{g\in\MCG/\MCH}\mu\left(\tilde{g}\MCH\right)g,\]
where for any $g\in\MCG/\MCH$, $\tilde{g}\in\MCG$ is an arbitrary preimage of $g$.
Then we have a canonical isomorphism
\[\MCM\left(\MCG,A\right) \xymatrix{\ar[r]^-{\sim} &} A\left[\left[\MCG\right]\right], \quad \mu \xymatrix{\ar@{|->}[r] &} \underline{\mu}, \]
and for any $\mu\in\MCM\left(\MCG,A\right)$ and any $\GGs\in\MCG$, we have
$\underline{\GGs_\ast\mu} = \GGs\underline{\mu}$.
Also we mention that if $\tilde{\pi}:A[[\MCG_1]] \rightarrow A[[\MCG_2]]$ is the canonical morphism defined by $\pi$, then we have the following commutative squares,
\[\xymatrix{
\MCM\left(\MCG,A\right) \ar[r]^-{\sim} \ar[d]_{\pi_\ast} & A\left[\left[\MCG\right]\right] \ar[d]^{\tilde{\pi}} \\
\MCM\left(\MCG',A\right) \ar[r]^-{\sim} & A\left[\left[\MCG'\right]\right] 
} \quad\text{and}\quad \xymatrix{
\MCM\left(\MCG,A\right) \ar[r]^-{\sim} & A\left[\left[\MCG\right]\right] & \Sigma h\\
\MCM\left(\MCG',A\right) \ar[r]^-{\sim} \ar[u]^{\pi^\sharp} & A\left[\left[\MCG'\right]\right] \ar[u] & g \ar@{|->}[u] 
}\]
where for all $g\in\MCG'$, $\Sigma h$ is the sum over all $h\in\MCG$ such that $\pi(h)=g$.

\begin{pr}\label{impibullet}
Let $\pi:\MCG_1\rightarrow\MCG_2$ be an open morphism of profinite groups, such that $\Ker(\pi)$ is finite.
The morphism $\pi^\sharp:\MCM\left(\MCG_2,A\right)\rightarrow\MCM\left(\MCG_1,A\right)$ is injective if and only if $\pi$ is surjective.
Moreover if $\#\left(\Ker(\pi)\right)$ is not a zero divisor in $A$, then the image of $\pi^\sharp$ is $\MCM\left(\MCG_1,A\right)^{\Ker(\pi)}$.
\end{pr}

\noindent\textsl{Proof.}
Let $\mu_2\in\MCM\left(\MCG_2,A\right)$.
From (\ref{defpibullet}) it is straightforward to check that $\pi^\sharp\mu_2=0$ if and only if $\mu_2(X)=0$ for all $X\in\MFX(\IM(\pi))$, and then we deduce that $\pi^\sharp$ is injective if and only if $\pi$ is surjective.
For any $\GGs\in\Ker(\pi)$, any $g\in\MCG_1$, and any open subgroup $\MCH$ of $\MCG_1$, we have
\[\begin{array} {cclcc}
\GGs_\ast\pi^\sharp\mu_2 \left(g\MCH\right) & = & \#\left(\MCH\cap\Ker(\pi)\right) \mu_2 \left( \pi\left(\GGs^{-1}g\MCH\right) \right) & & \\
 & = & \#\left(\MCH\cap\Ker(\pi)\right) \mu_2 \left( \pi\left(g\MCH\right) \right) & = & \pi^\sharp\mu_2 \left(g\MCH\right),
\end{array}\]
hence $\GGs_\ast\pi^\sharp\mu_2 = \pi^\sharp\mu_2$, and $\IM\left(\pi^\sharp\right) \subseteq \MCM\left(\MCG_1,A\right)^{\Ker(\pi)}$.

Now let $\mu_1\in\MCM\left(\MCG_1,A\right)^{\Ker(\pi)}$.
Let $\MCH$ be an open subgroup of $\IM(\pi)$, and $g\in\MCG_1$.
Let $\MCW$ be an open normal subgroup of $\pi^{-1}(\MCH)$ such that $\MCW\cap\Ker(\pi)$ is trivial.
Let $R$ be a complete representative system of $\pi^{-1}(\MCH)$ modulo $\MCW\Ker(\pi)$.
Then $(\GGs r)_{(\GGs,r)\in\Ker(\pi)\times R}$ is a complete representative system of $\pi^{-1}(\MCH)$ modulo $\MCW$, and we have
\[\begin{array}{cclcc}
\mu_1\left(\pi^{-1}\left(\MCH\right)g\right) & = & \sum_{(\GGs,r)\in\Ker(\pi)\times R} \mu_1\left(\GGs r\MCW g\right) & = & \sum_{(\GGs,r)\in\Ker(\pi)\times R} \left(\GGs^{-1}\right)_\ast \mu_1 \left( r\MCW g \right) \\
 & = & \sum_{(\GGs,r)\in\Ker(\pi)\times R} \mu_1\left(r\MCW g\right) & = & \#\left(\Ker(\pi)\right) \sum_{r\in R} \mu_1\left(r\MCW g\right).
\end{array}\]
Hence $\pi_\ast\mu_1$ takes values in $\#\left(\Ker(\pi)\right)A$, and we deduce $\mu_1 = \pi^\sharp \left(\#\left(\Ker(\pi)\right)^{-1}\pi_\ast\mu_1\right)$ from (\ref{pibulletpiast}).
\hfill $\square$
\\

Now assume $A:=\MCO_L$ for some complete subfield $L$ of $\MBC_p$, finitely ramified over $\MBQ_p$.
An $A$-distribution on $\MCG$ is called a measure.
Let $\mu\in\MCM\left(\MCG,A\right)$ be such a measure, and let $V$ be a complete topological $A$-module, such that the open submodules of $V$ form a neighborhood basis for $0$.
Let $\MSC\left(\MCG,V\right)$ be the $A$-module of continuous maps from $\MCG$ to $V$, equipped with the uniform convergence topology.
For any $X\in\MFX(\MCG)$, we denote by $1_X:\MCG\rightarrow A$ the map such that $1_X(x)=1$ for $x\in X$ and $1_X(x)=0$ for $x\in \MCG\setminus X$.
Then there is a unique continuous $A$-linear map 
\[\xymatrix{\MSC\left(\MCG,V\right)\ar[r] & V, & f \ar@{|->}[r] &} \int f(t).\dd\mu(t),\]
such that for all $X\in\MFX(\MCG)$ and all $v\in V$, $\int 1_X(t)v.\dd\mu(t) = \mu(X)v$.
For $X\in\MFX(\MCG)$ and $f\in\MSC\left(\MCG,V\right)$, we write $\int_Xf.\dd\mu$ for $\int1_Xf.\dd\mu$.
Then for $\GGs\in\MCG$, we have
\begin{equation}
\label{GGsheufg}
\int_Xf(t).\dd\GGs_\ast\mu(t) \; = \; \int_{\GGs^{-1}X}f(\GGs t).\dd\mu(t).
\end{equation}
For any $\GGa\in\MBZ_p$, recall that $(1+T)^\GGa:=1+\sum_{n=1}^\infty \frac {\GGa (\GGa-1) ... (\GGa+1-n)} {n!} T^n$ belongs to $\MBZ_p[[T]]$.
Then for $\mu\in\MCM\left(\GGG,A\right)$, we have
\begin{equation}
\label{socnfyhhf}
\underline{\mu} = \int \left(1+T\right)^{\GGk(\GGs)} .\dd\mu(\GGs) \quad \text{in} \quad A[[T]],
\end{equation}
where $\GGk:\GGG\rightarrow\MBZ_p$ is the unique isomorphism of profinite groups such that $\GGk(\GGg)=1$.
Moreover, if we write $\MFm_{\MBC_p}$ for the maximal ideal of $\MCO_{\MBC_p}$, then for any $x\in\MFm_{\MBC_p}$ we have
\begin{equation}
\label{socnfyhhfdeux}
\underline{\mu}(x) = \int \left(1+x\right)^{\GGk(\GGs)} .\dd\mu(\GGs) \quad \text{in} \quad \MBC_p.
\end{equation}

\section{Elliptic units.}\label{sectionellipticunits}

For $L$ and $L'$ two $\MBZ$-lattices of $\MBC$ such that $L\subseteq L'$ and $[L':L]$ is prime to $6$, we denote by $z\mapsto\psi\left(z;L,L'\right)$ the elliptic function defined in \cite{robert90}.
For $\MFm$ a nonzero proper ideal of $\MCO_k$, and $\MFa$ a nonzero ideal of $\MCO_k$ prime to $6\MFm$, G.\,Robert proved that $\psi\left(1;\MFm,\MFa^{-1}\MFm\right) \in k(\MFm)$, where $k(\MFm)$ is the ray class field of $k$, modulo $\MFm$.
Let $S(\MFm)$ be the set of maximal ideals of $\MCO_k$ which divide $\MFm$.
Then $\psi\left(1;\MFm,\MFa^{-1}\MFm\right)$ is a unit if and only if $\vert S(\MFm)\vert = 1$.
More precisely, if we denote by $w_\MFm$ the number of roots of unity of $k$ which are congruent to $1$ modulo $\MFm$, and if we write $w_k$ for the number of roots of unity in $k$, then by \cite[Corollaire 1.3, (iv)]{robert89}, we have
\begin{equation}
\label{unitellipticisunit}
\psi\left(1;\MFm,\MFa^{-1}\MFm\right) \MCO_{k(\MFm)} = \lA \begin{array}{lll} (1) & \text{if} & 2\leq\vert S(\MFm)\vert \\ (\MFq)_{k(\MFm)}^{w_\MFm\left(N(\MFa)-1\right)/w_k} & \text{if} & S(\MFm) = \{\MFq\}, \end{array} \right.
\end{equation}
where $N(\MFa):=\#\left(\MCO_k/\MFa\right)$, and where $(\MFq)_{k(\MFm)}$ is the product of the prime ideals of $\MCO_{k(\MFm)}$ which lie above $\MFq$.
Moreover, if $\GVf_\MFm(1)$ is the Robert-Ramachandra invariant, as defined in \cite[p15]{robert73}, or in \cite[p55]{deshalit87}, we have by \cite[Corollaire 1.3, (iii)]{robert89}
\begin{equation}
\label{ellipticunitpsiandphi}
\psi\left(1;\MFm,\MFa^{-1}\MFm\right) ^{12m} = \GVf_\MFm(1)^{N(\MFa)-\left(\MFa,k(\MFm)/k\right)},
\end{equation}
where $m$ is the positive generator of $\MFm\cap\MBZ$, and $\left(\MFa,k(\MFm)/k\right)$ is the Artin automorphism of $k(\MFm)/k$ defined by $\MFa$.
If $\MFa$ is prime to $6\MFm\MFq$, then by \cite[Corollaire 1.3, (ii-1)]{robert89} we have
\begin{equation}
\label{ellipticunit}
N_{k(\MFm\MFq)/k(\MFm)} 
\left( \psi\left(1;\MFm\MFq,\MFa^{-1}\MFm\MFq\right) \right) ^{w_\MFm w_{\MFm\MFq}^{-1}} =
\lA
\begin{array}{lcl}
\psi\left(1;\MFm,\MFa^{-1}\MFm\right)^{1-\left(\MFq,k(\MFm)/k\right)^{-1}} & \mathrm{if} & \MFq\nmid\MFm,\\
\psi\left(1;\MFm,\MFa^{-1}\MFm\right) \vphantom{\left(I^J\right)^K} & \mathrm{if} & \MFq\mid\MFm.
\end{array}
\right.
\end{equation}

\begin{df}
Let $F\subset\MBC$ be a finite abelian extension of $k$, and write $\mu(F)$ for the group of roots of unity in $F$.
Let $\MFm$ be a nonzero proper ideal of $\MCO_k$.
We define the $\MBZ\left[\Gal(F/k)\right]$-submodule $\Psi(F,\MFm)$ of $F^\times$, generated by the $w_\MFm$-roots of all $N_{k(\MFm)/k(\MFm)\cap F} \left( \psi\left(1;\MFm,\MFa^{-1}\MFm\right) \right)$, where $\MFa$ is any nonzero ideal of $\MCO_k$ prime to $6\MFm$.
Also, we set $\Psi'(F,\MFm):=\MCO_F^\times\cap\Psi(F,\MFm)$.

Then, we let $C_F$ be the group generated by $\mu(F)$ and by all $\Psi'(F,\MFm)$, for any nonzero proper ideal $\MFm$ of $\MCO_k$.
\end{df}

\begin{rmq}\label{astrototoj}
Let $\MFm$ and $\MFg$ be two nonzero proper ideals of $\MCO_k$, such that the conductor of $F/k$ divides $\MFm$.
If $\MFg\wedge\MFm=1$, then $\Psi'\left(F,\MFg\right) \subseteq C_F\cap\MCO_{k(1)}^\times$.
Else by (\ref{ellipticunit}) we have $\Psi'\left(F,\MFg\right) \subseteq \Psi'\left(F,\MFg\wedge\MFm\right)$.
\end{rmq}

We define $\MCC_n:=\MBZ_p\otimes_\MBZ C_n$, and $\MCC_\infty := \varprojlim \left(\MCC_n\right)$, projective limit under the norm maps.
For any nonzero ideal $\MFg$ of $\MCO_k$, we define 
\[\Psi\left(K_n,\MFg\MFp^\infty\right) := \mathop{\cup}_{i=1}^\infty \Psi\left(K_n,\MFg\MFp^i\right) \quad 
\text{and} \quad \Psi'\left(K_n,\MFg\MFp^\infty\right) := \mathop{\cup}_{i=1}^\infty \Psi'\left(K_n,\MFg\MFp^i\right).\]
Then the projective limits under the norm maps are denoted by
\[\overline{\Psi}\left(K_\infty,\MFg\MFp^\infty\right) := \varprojlim \left( \MBZ_p\otimes_\MBZ\Psi\left(K_n,\MFg\MFp^\infty\right) \right), \quad 
\overline{\Psi}'\left(K_\infty,\MFg\MFp^\infty\right) := \varprojlim \left( \MBZ_p\otimes_\MBZ\Psi'\left(K_n,\MFg\MFp^\infty\right) \right).\]

Let us write $\MCI$ for the set of nonzero ideals of $\MCO_k$ which are prime to $\MFp$.
For $\MFg\in\MCI$, we set $K_{\MFg,\infty} := k\left(\MFg\MFp^\infty\right)=\mathop{\cup}_{n\in\MBN}k\left(\MFg\MFp^n\right)$, and $G_{\MFg,\infty}:=\Gal\left(K_{\MFg,\infty}/k\right)$.
Then we write $\GGG_\MFg$ for the torsion subgroup of $G_{\MFg,\infty}$.
We denote by $\MCI'$ the subset of $\MCI$ containing all the $\MFg\in\MCI$ such that $w_\MFg=1$.
In the sequel, we fix once and for all $\MFf\in\MCI'$ such that $K_\infty\subseteq K_{\MFf,\infty}$.
We choose arbitrarily a subgroup of $G_{\MFf,\infty}$, isomorphic to $\MBZ_p$, such that its image in $G_\infty$ is $\GGG$.
Then for any $\MFg\in\MCI$ such that $\MFg|\MFf$, we have the decomposition $G_{\MFg,\infty}=G_\MFg\times\GGG$.

\begin{rmq}\label{Dragoon}
From Remark \ref{astrototoj}, $\MCC_\infty$ is generated by all the $\overline{\Psi}'\left(K_\infty,\MFg\MFp^\infty\right)$, where $\MFg\in\MCI$ is such that $\MFg\vert\MFf$.
\end{rmq}

From (\ref{ellipticunit}), for $\MFg\in\MCI$ such that $\MFg\vert\MFf$, and for any nonzero ideal $\MFa$ of $\MCO_k$ which is prime to $6\MFg\MFp$, 
there is a unique 
\[\psi\la\MFg,\MFa\ra \in \overline{\Psi}\left(K_{\MFg,\infty},\MFg\MFp^\infty\right)\] 
such that for large enough $n\in\MBN$, 
the canonical image of $\psi\la\MFg,\MFa\ra$ in $\MBZ_p\otimes_\MBZ\Psi\left(k\left(\MFg\MFp^n\right),\MFg\MFp^\infty\right)$ is $1\otimes\psi\left(1;\MFg\MFp^n,\MFa^{-1}\MFg\MFp^n\right)$.

\section{From semilocal units to measures.}

Let $\MBQ_p^{\nr}\subseteq\MBC_p$ be the maximal unramified algebraic extension of $\MBQ_p$, and let $L$ be the completion of $\MBQ_p^{\nr}$.
We denote by $\MCO_\MFf$ the ring $\MCO_L\left[\zeta\right]$, where $\zeta$ is any primitive $\left[K_{\MFf,0} : k\right]$-th root of unity in $\MBC_p$.
For all $\left(\MFg_1,\MFg_2\right)\in\MCI^2$ such that $\MFg_1\vert\MFg_2$, we denote by $\pi_{\MFg_2,\MFg_1} : G_{\MFg_2,\infty} \rightarrow G_{\MFg_1,\infty}$ the restriction map.
We write $N_{\MFg_2,\MFg_1} : \MCO_\MFf \widehat{\otimes}_{\MBZ_p} \MCU_{\MFg_2,\infty} \rightarrow \MCO_\MFf \widehat{\otimes}_{\MBZ_p} \MCU_{\MFg_1,\infty}$ for the norm map and we write $\GGu_{\MFg_2,\MFg_1} : \MCO_\MFf \widehat{\otimes}_{\MBZ_p} \MCU_{\MFg_1,\infty} \rightarrow \MCO_\MFf \widehat{\otimes}_{\MBZ_p} \MCU_{\MFg_2,\infty}$ for the canonical injection.

For all $\MFg\in\MCI'$, de Shalit defined in \cite[I.3.4, II.4.6, and II.4.7]{deshalit87} an injective morphism of $\MBZ_p\left[\left[G_\infty\right]\right]$-modules $i_\MFg^0 : \MCU_{\MFg,\infty} \rightarrow \MCM\left( G_{\MFg,\infty} , \MCO_L \right)$, which we extend by linearity into a morphism $i_\MFg$ from $\MCO_\MFf \widehat{\otimes}_{\MBZ_p} \MCU_{\MFg,\infty}$ to $\MCM\left( G_{\MFg,\infty} , \MCO_\MFf \right)$.

\begin{lm}\label{sixcentquin}
There is a unique way to extend the family $\left(i_\MFg\right)_{\MFg\in\MCI'}$ to $\MCI$ such that for all $\left(\MFg_1,\MFg_2\right)\in\MCI^2$, the following squares are commutative,
\begin{equation}
\label{relationiuhfreg}
\xymatrix{\MCO_\MFf \widehat{\otimes}_{\MBZ_p} \MCU_{\MFg_2,\infty} \ar[d]^-{N_{\MFg_2,\MFg_1}} \ar[r]^-{i_{\MFg_2}} & \MCM\left(G_{\MFg_2,\infty},\MCO_\MFf\right) \ar[d]^-{ \left(\pi_{\MFg_2,\MFg_1}\right)_\ast } \\
\MCO_\MFf \widehat{\otimes}_{\MBZ_p} \MCU_{\MFg_1,\infty} \ar[r]^-{i_{\MFg_1}} & \MCM\left(G_{\MFg_1,\infty},\MCO_\MFf\right)}
\quad \text{and} \quad
\xymatrix{\MCO_\MFf \widehat{\otimes}_{\MBZ_p} \MCU_{\MFg_1,\infty} \ar[d]^-{\GGu_{\MFg_2,\MFg_1}} \ar[r]^-{i_{\MFg_1}} & \MCM\left(G_{\MFg_1,\infty},\MCO_\MFf\right) \ar[d]^-{ \left(\pi_{\MFg_2,\MFg_1}\right)^\sharp } \\
\MCO_\MFf \widehat{\otimes}_{\MBZ_p} \MCU_{\MFg_2,\infty} \ar[r]^-{i_{\MFg_2}} & \MCM\left(G_{\MFg_2,\infty},\MCO_\MFf\right)}
\end{equation}
\end{lm}

\noindent\textsl{Proof.}
This was proved by de Shalit in the case $p\neq2$ (see \cite[III.1.2 and III.1.3]{deshalit87}).
Let $\MFg_1\in\MCI\setminus\MCI'$, and let $\MFg_2\in\MCI'$ be such that $\MFg_1\vert\MFg_2$.
When $p\neq2$ de Shalit uses the surjectivity of $N_{\MFg_2,\MFg_1}$ in order to construct $i_{\MFg_1}$.
If $p=2$, $N_{\MFg_2,\MFg_1}$ may not be surjective.
However we have 
$\IM\left(i_{\MFg_2}\circ\GGu_{\MFg_2,\MFg_1}\right) \subseteq \MCM \left( G_{\MFg_2,\infty},\MCO_\MFf \right) ^{\Ker\left(\pi_{\MFg_2,\MFg_1}\right)}$.
But by Proposition \ref{impibullet}, $\left(\pi_{\MFg_2,\MFg_1}\right)^\sharp$ is injective and $\IM\left(\pi_{\MFg_2,\MFg_1}\right)^\sharp = \MCM \left( G_{\MFg_2,\infty},\MCO_\MFf \right) ^{\Ker\left(\pi_{\MFg_2,\MFg_1}\right)}$.
Hence there is a unique map $i_{\MFg_1}$ such that the right hand square of (\ref{relationiuhfreg}) is commutative.
The rest of the proof is identical to \cite{deshalit87}.
\hfill $\square$
\\

\begin{lm}
For all $\MFg\in\MCI$, $i_\MFg$ is an injective pseudo-isomorphism of $\MCO_\MFf[[T]]$-modules.
\end{lm}

\noindent\textsl{Proof.}
Let $\MFg_1\in\MCI$, and let $\MFg_2\in\MCI'$ be such that $\MFg_1\vert\MFg_2$.
Then $\left(\pi_{\MFg_2,\MFg_1}\right)^\sharp$, $\GGu_{\MFg_2,\MFg_1}$, and $i_{\MFg_2}$ are injective, and by (\ref{relationiuhfreg}) we deduce the injectivity of $i_{\MFg_1}$.

By class field theory, one can show that for any prime $\MFq$ of $K_{\MFg_2,\infty}$ above $\MFp$, the number of $p$-power roots of unity in $\left(K_{\MFg_2,n}\right)_\MFq$ is bounded independantly of $n$ (see \cite[Lemma 2.1]{viguie11b}).
Then it follows from \cite[I.3.7, Theorem]{deshalit87} that $i_{\MFg_2}$ is a pseudo-isomorphism.
Since $\left( \MCO_\MFf \widehat{\otimes}_{\MBZ_p} \MCU_{\MFg_2,\infty} \right) ^ {\Ker\left(\pi_{\MFg_2,\MFg_1}\right)} = \MCO_\MFf \widehat{\otimes}_{\MBZ_p} \MCU_{\MFg_1,\infty}$ and since $\left(\pi_{\MFg_2,\MFg_1}\right)^\sharp$ is injective, 
it follows from (\ref{relationiuhfreg}) that $\MCM\left(G_{\MFg_1,\infty} , \MCO_\MFf\right) / \IM\left(i_{\MFg_1}\right)$ is a submodule of $\MCM\left(G_{\MFg_2,\infty} , \MCO_\MFf\right) / \IM\left(i_{\MFg_2}\right)$, which is pseudo-nul since $i_{\MFg_2}$ is a pseudo-isomorphism.
\hfill $\square$
\\

An element of the total fraction ring of $\MCM\left( G_{\MFg,\infty} , \MCO_L\right)$ is called an $\MCO_L$-pseudo-measure.
For $\MFg\in\MCI$, let $\mu(\MFg)$ be the $\MCO_L$-pseudo-measure on $G_{\MFg,\infty}$ defined in \cite[II.4.12, Theorem]{deshalit87}.
It is a measure if $\MFg\neq(1)$, and $\GGa\mu(1)$ is a measure for all $\GGa\in\MCJ_{(1)}$, where we write $\MCJ_{(1)}$ for the augmentation ideal of $\MCO_\MFf\left[\left[G_{(1),\infty}\right]\right]$.
By definition of $\mu(\MFg)$, we have
\begin{equation}
\label{igfhgffj}
i_\MFg\left(\psi\la\MFg,\MFa\ra\right) = \left(\left(\MFa,K_{\MFg,\infty}/k\right) - N(\MFa)\right) \mu(\MFg).
\end{equation}
Moreover, for $\left(\MFg_1,\MFg_2\right)\in\MCI^2$ such that $\MFg_1\vert\MFg_2$, we have
\begin{equation}
\label{xomftnf}
\left(\pi_{\MFg_2,\MFg_1}\right)_\ast\mu\left(\MFg_2\right) = \prod_{ \substack{ \MFl \, \text{prime of} \, \MCO_k \\ \MFl\vert\MFg_2 \, \text{and} \, \MFl\nmid\MFg_1}} \left(1-\left(\MFl,K_{\MFg_1,\infty}/k\right)^{-1}\right)\mu\left(\MFg_1\right).
\end{equation}

\begin{lm}\label{NIUbhgREeeE}
For $\MFg\in\MCI$, we denote by $\mu_{p^\infty} \left(K_{\MFg,\infty}\right)$ the group of $p$-power roots of unity in $K_{\MFg,\infty}$.
Then we have
\[i_\MFg \left( \MCO_\MFf \widehat{\otimes}_{\MBZ_p} \overline{\Psi}' \left( K_{\MFg,\infty} , \MFg\MFp^\infty \right) \right) = \MCJ_\MFg \mu(\MFg),\]
where $\MCJ_\MFg$ is the annihilator of the $\MCO_\MFf\left[\left[G_{\MFg,\infty}\right]\right]$-module $\MCO_\MFf \widehat{\otimes}_{\MBZ_p} \mu_{p^\infty} \left(K_{\MFg,\infty}\right)$ if $\MFg\neq(1)$, 
and where $\MCJ_{(1)}$ is the augmentation ideal of $\MCO_\MFf\left[\left[G_{(1),\infty}\right]\right]$.
\end{lm}

\noindent\textsl{Proof.}
We refer the reader to \cite[III.1.4]{deshalit87}.
\hfill $\square$
\\

\section{Generation of the characteristic ideal.}

For any $\MFg\in\MCI$ such that $\MFg\vert\MFf$, and any irreducible ($\MBC$ or $\MBC_p$) character $\chi$ of $G_\MFg$, let $\MFf_\chi\in\MCI$ be such that the conductor of $\chi$ is $\MFf_\chi\MFp^n$ for some $n\in\MBN$.
Then $\chi$ defines a character on $G_{\MFf_\chi}$, which we denote by $\chi_0$.
We have
\[\MCO_\MFf\left[\left[G_{\MFg,\infty}\right]\right]_\chi \simeq \MCO_\MFf\left[\left[\GGG\right]\right] \quad\text{and}\quad \MCM\left(G_{\MFg,\infty},\MCO_\MFf\right)_\chi \simeq \MCM\left(\GGG,\MCO_\MFf\right),\]
where the isomorphisms are induced by the following maps,
\[\tilde{\chi} : \MCO_\MFf\left[\left[G_{\MFg,\infty}\right]\right] \rightarrow \MCO_\MFf\left[\left[\GGG\right]\right] \quad\text{and}\quad \chi' : \MCM\left(G_{\MFg,\infty},\MCO_\MFf\right) \rightarrow \MCM\left(\GGG,\MCO_\MFf\right),\]
such that for any $(g,\GGs)\in G_{\MFg}\times\GGG$, $\tilde{\chi}(\GGs g)=\chi(g)\GGs$, and such that for any $\mu\in\MCM\left(G_{\MFg,\infty},\MCO_\MFf\right)$, $\underline{\chi'(\mu)} = \tilde{\chi}\left(\underline{\mu}\right)$.
Moreover, remark that we have
\begin{eqnarray}
\label{equationdjvjfnfjf}
& & \chi'(\mu) \; = \; \chi'_0\left( \left( \pi_{\MFg,\MFf_\chi} \right)_\ast \mu\right) \quad \text{for all} \quad \mu\in\MCM\left(G_{\MFg,\infty},\MCO_\MFf\right), \\
 & \text{and} & \nonumber\\
\label{efivnggy}
& & \chi'\circ\left(\pi_{\MFg,\MFh}\right)^\sharp \; = \; 0 \quad \text{for all $\MFh\in\MCI$ such that $\MFh\neq\MFf_\chi$ and $\MFh\vert\MFf_\chi$.}
\end{eqnarray}

For any finite group $\MCG$, any irreducible $\MBC_p$-character $\chi$ of $\MCG$, and any morphism $f:M\rightarrow N$ of $\MCO_\MFf \left[\MCG\right]$-modules, we denote by $f_\chi : M_\chi\rightarrow N_\chi$ the morphism defined by $f$.
For any $x\in M$, we write $x_\chi$ for the canonical image of $x$ in $M_\chi$.

\begin{lm}\label{kdugng}
Let $\MFg\in\MCI$ be such that $\MFg\vert\MFf$.
Let $\chi\neq1$ be an irreducible $\MBC_p$-character of $G_\MFg$.
Then
\[\left(i_{\MFg}\right)_\chi \left( \left(\MCO_\MFf \widehat{\otimes}_{\MBZ_p} \MCC_{\MFg,\infty}\right)_\chi \right) \subseteq \left(i_{\MFf_\chi}\right)_{\chi_0} \left( \left(\MCO_\MFf \widehat{\otimes}_{\MBZ_p} \overline{\Psi}'\left( K_{\MFf_\chi,\infty} , \MFf_\chi\MFp^\infty \right) \right)_{\chi_0} \right),\]
and the quotient is a pseudo-nul $\MCO_\MFf[[T]]$-module.
\end{lm}

\noindent\textsl{Proof.}
Let $\MFh\in\MCI$ be such that $\MFh\vert\MFg$, and let $x\in\overline{\Psi}'\left( K_{\MFg,\infty} , \MFh\MFp^\infty \right)$.
From Remark \ref{astrototoj}, there is $y\in\overline{\Psi}'\left( K_{\MFh\wedge\MFf_\chi,\infty} , \left( \MFh\wedge\MFf_\chi \right)\MFp^\infty \right)$ such that $N_{\MFg,\MFf_\chi}(x) = \GGu_{\MFf_\chi,\MFh\wedge\MFf_\chi}(y)$.
From  (\ref{equationdjvjfnfjf}), and then from (\ref{relationiuhfreg}), one has
\begin{eqnarray}
\label{equassdftion}
\left(i_{\MFg}\right)_\chi \left(x_\chi\right) = \chi' \circ i_{\MFg}(x) = \chi'_0 \circ \left(\pi_{\MFg,\MFf_\chi}\right)_\ast \circ i_{\MFg}(x) & = & \chi'_0 \circ i_{\MFf_\chi} \circ N_{\MFg,\MFf_\chi}(x) = \chi'_0 \circ i_{\MFf_\chi} \circ \GGu_{\MFf_\chi,\MFh\wedge\MFf_\chi}(y)\nonumber\\
 & = &\chi'_0 \circ \left(\pi_{\MFf_\chi,\MFh\wedge\MFf_\chi}\right)^\sharp \circ i_{\MFh\wedge\MFf_\chi}(y).
\end{eqnarray}
From (\ref{efivnggy}) and (\ref{equassdftion}), we deduce $\left(i_{\MFg}\right)_\chi \left(x_\chi\right) = 0$ if $\MFf_\chi\nmid\MFh$, and $\left(i_{\MFg}\right)_\chi \left(x_\chi\right) = \chi'_0 \circ i_{\MFf_\chi}(y) = \left(i_{\MFf_\chi}\right)_{\chi_0} \left(y_{\chi_0}\right)$ if $\MFf_\chi\vert\MFh$.
By Remark \ref{Dragoon}, this states the inclusion $\MCB\subseteq\MCA$, where we set
\[\MCA := \left(i_{\MFf_\chi}\right)_{\chi_0} \left( \left(\MCO_\MFf \widehat{\otimes}_{\MBZ_p} \overline{\Psi}'\left( K_{\MFf_\chi,\infty} , \MFf_\chi\MFp^\infty \right) \right)_{\chi_0} \right) \quad\text{and}\quad \MCB := \left(i_{\MFg}\right)_\chi \left( \left(\MCO_\MFf \widehat{\otimes}_{\MBZ_p} \MCC_{\MFg,\infty}\right)_\chi \right).\]
Let $m:=\left[k\left(\MFg\MFp^\infty\right) : k\left(\MFf_\chi\MFp^\infty\right)\right]$, and let $x\in \overline{\Psi}' \left(K_{\MFf_\chi,\infty} , \MFf_\chi\MFp^\infty\right)$.
Then $mx = N_{\MFg,\MFf_\chi} \circ \GGu_{\MFg,\MFf_\chi}(x)$, and from (\ref{relationiuhfreg}) and (\ref{equationdjvjfnfjf}), we obtain
\begin{eqnarray*}
m\left(i_{\MFf_\chi}\right)_{\chi_0} \left(x_{\chi_0}\right) = \chi'_0\circ i_{\MFf_\chi} \circ N_{\MFg,\MFf_\chi} \circ \GGu_{\MFg,\MFf_\chi}(x) & = & \chi'_0 \circ \left(\pi_{\MFg,\MFf_\chi}\right)_\ast \circ i_\MFg \circ \GGu_{\MFg,\MFf_\chi}(x) \\
& = & \chi' \circ i_\MFg \circ \GGu_{\MFg,\MFf_\chi}(x) = \left(i_\MFg\right)_\chi \left( \GGu_{\MFg,\MFf_\chi}(x)_\chi \right), 
\end{eqnarray*}
and we deduce that $m$ annihilates $\MCA/\MCB$.
Let $\GGa := \prod_{ \substack{ \MFl \, \text{prime of} \, \MCO_k \\ \MFl\vert\MFg \, \text{and} \, \MFl\nmid\MFf_\chi}} \left(1-\tilde{\chi}_0\left(\GGs_\MFl^{-1}\right)\right)$, where $\GGs_\MFl$ is the Fr\"obenius of $\MFl$ in $K_{\MFf_\chi,\infty}/k$.
Let $x\in\overline{\Psi}' \left(K_{\MFf_\chi,\infty} , \MFf_\chi\MFp^\infty\right)$.
From (\ref{ellipticunit}), there is $y\in\overline{\Psi}' \left(K_{\MFg,\infty} , \MFg\MFp^\infty\right)$ such that $\GGa x = N_{\MFg,\MFf_\chi}(y)$.
Then by (\ref{relationiuhfreg}) and (\ref{equationdjvjfnfjf}), we have
\begin{eqnarray*}
\GGa \left(i_{\MFf_\chi}\right)_{\chi_0} \left(x_{\chi_0}\right) = \chi'_0 \circ i_{\MFf_\chi} \circ N_{\MFg,\MFf_\chi} (y) = \chi'_0 \circ \left(\pi_{\MFg,\MFf_\chi}\right)_\ast \circ i_\MFg (y) & = & \chi' \circ i_\MFg (y) \\
 & = & \left(i_\MFg\right)_\chi \left(y_\chi\right).
\end{eqnarray*}
Hence $\GGa$ annihilates $\MCA/\MCB$.
As a particular case, if there is no maximal ideal $\MFl$ of $\MCO_k$ such that $\MFl\vert\MFg$ and $\MFl\nmid\MFf_\chi$, then $\GGa=1$, $\MCA=\MCB$, and Lemma \ref{kdugng} is proved in this case.
Now assume that there is a maximal ideal $\MFl$ of $\MCO_k$ such that $\MFl\vert\MFg$ and $\MFl\nmid\MFf_\chi$.
By class field theory, the decomposition group of $\MFl$ in $K_{\MFf_\chi,\infty}/k$ has a finite index in $\Gal\left( K_{\MFf_\chi,\infty} / k \right)$.
Hence $\GGs_\MFl\notin G_{\MFf_\chi}$, and there are a topological generator $\tilde{\GGg}$ of $\GGG$, $n\in\MBN$, and $g\in G_{\MFf_\chi}$ such that $\GGs_\MFl^{-1}=g\tilde{\GGg}^{p^n}$.
Then 
\begin{equation}
\label{dujdjkf}
1-\tilde{\chi}_0\left(\GGs_\MFl^{-1}\right) = 1-\chi_0(g) \tilde{\GGg}^{p^n} = 1 - \chi_0(g) \sum_{i=0}^{p^n} \left( \begin{array}{c} p^n \\ i \end{array} \right) \tilde{T}^j, 
\end{equation}
where $\tilde{T} := \tilde{\GGg}-1$.
Since $m$ and $\chi_0(g)$ are coprime, and since $-\chi_0(g)$ is the coefficient of $\tilde{T}^{p^n}$ in the decomposition (\ref{dujdjkf}), we deduce that $m$ and $1-\tilde{\chi}_0\left(\GGs_\MFl^{-1}\right)$ are coprime.
Then $m$ and $\GGa$ are coprime, and annihilates $\MCA/\MCB$, so that Lemma \ref{kdugng} follows.
\hfill $\square$
\\

\begin{lm}\label{cccodglh}
Let $\MFg\in\MCI$ be such that $\MFg\vert\MFf$.
Let $\chi\neq1$ be an irreducible $\MBC_p$-character of $G_\MFg$.

$\mathrm{(i)}$ If $p\neq2$ or if $w_\MFg = w_{\MFf_\chi}$, then $\IM\left(i_\MFg\right)_\chi = \IM\left(i_{\MFf_\chi}\right)_{\chi_0}$.

$\mathrm{(ii)}$ If $p=2$, then $\IM\left(i_\MFg\right)_\chi \subseteq \IM\left(i_{\MFf_\chi}\right)_{\chi_0}$, and the quotient is annihilated by $2$.
\end{lm}

\noindent\textsl{Proof.}
For $x\in\MCU_{\MFg,\infty}$, by (\ref{equationdjvjfnfjf}) and (\ref{relationiuhfreg}), we have
\begin{equation}
\label{equationdhdddh}
\left(i_\MFg\right)_\chi\left(x_\chi\right) = \chi'_0 \circ \left(\pi_{\MFg,\MFf_\chi}\right)_\ast \circ i_\MFg (x) = \chi'_0 \circ i_{\MFf_\chi} \circ N_{\MFg,\MFf_\chi} (x) =  \left(i_{\MFf_\chi}\right)_{\chi_0} \left(N_{\MFg,\MFf_\chi} (x) _{\chi_0}\right). 
\end{equation}
We deduce $\IM\left(i_\MFg\right)_\chi \subseteq \IM\left(i_{\MFf_\chi}\right)_{\chi_0}$.
For $n$ large enough, the ramification index of the primes above $\MFp$ in $K_{\MFg,n}/K_{\MFf_\chi,n}$ is $w_{\MFf_\chi}w_{\MFg}^{-1}$.
If $p\neq 2$, then $w_{\MFf_\chi}w_{\MFg}^{-1}$ is prime to $p$.
Hence in case $\mathrm{(i)}$, $K_{\MFg,n}/K_{\MFf_\chi,n}$ is tamely ramified.
Then $N_{\MFg,\MFf_\chi}$ is a surjection from $\MCU_{\MFg,\infty}$ onto $\MCU_{\MFf_\chi,\infty}$, and we deduce $\IM\left(i_\MFg\right)_\chi \supseteq \IM\left(i_{\MFf_\chi}\right)_{\chi_0}$ from (\ref{equationdhdddh}).
If $p=2$, $\MCU_{\MFf_\chi,\infty} / N_{\MFg,\MFf_\chi} \left( \MCU_{\MFg,\infty} \right)$ is annihilated by $w_{\MFf_\chi}w_{\MFg}^{-1}$ which is $1$ or $2$, and we deduce $\mathrm{(ii)}$ from (\ref{equationdhdddh}).
\hfill $\square$
\\

For $p\neq2$, Theorems \ref{theodufnfkffl} and \ref{theodufnfkffldeux} below were already proved by de Shalit in \cite[III.1.10]{deshalit87}.

\begin{theo}\label{theodufnfkffl}
Let $\MFg\in\MCI$ be such that $\MFg\vert\MFf$.
Let $\MFu$ be a uniformizer of $\MCO_\MFf$.
Let $\chi\neq1$ be an irreducible $\MBC_p$-character of $G_\MFg$.

$\mathrm{(i)}$ If $p\neq2$ or if $w_\MFg = w_{\MFf_\chi}$, then $\Char _{\MCO_\MFf[[T]]} \left( \MCO_\MFf \widehat{\otimes}_{\MBZ_p} \left(\MCU_{\MFg,\infty} / \MCC_{\MFg,\infty}\right) \right)_\chi$ is generated by $\tilde{\chi}_0 \left( \underline{\mu\left(\MFf_\chi\right)} \right)$.

$\mathrm{(ii)}$ If $p=2$, then $\Char _{\MCO_\MFf[[T]]} \left( \MCO_\MFf \widehat{\otimes}_{\MBZ_p} \left(\MCU_{\MFg,\infty} / \MCC_{\MFg,\infty}\right) \right)_\chi$ is generated by $\MFu^{-m_\chi}\tilde{\chi}_0 \left( \underline{\mu\left(\MFf_\chi\right)} \right)$, for some $m_\chi\in\MBN$.

(In case $\MFf_\chi=(1)$, we have expanded $\tilde{\chi}_0$ to the total fraction ring of $\MCO_\MFf\left[\left[G_{(1),\infty}\right]\right]$ and to the fraction field of $\MCO_\MFf[[\GGG]]$. 
We still have $\tilde{\chi}_0 \left( \underline{\mu\left(1\right)} \right)\in\MCO_\MFf[[\GGG]]$.)
\end{theo}

\noindent\textsl{Proof.}
Let us set $\tilde{\MCC}_\MFg := \left(\MCO_\MFf\widehat{\otimes}_{\MBZ_p}\MCC_{\MFg,\infty}\right)_\chi$.
We have the tautological exact sequence below,
\begin{equation}
\label{ddfff}
\xymatrix{0 \ar[r] & \IM \left(i_{\MFg}\right)_\chi / \left(i_{\MFg}\right)_\chi\left(\tilde{\MCC}_\MFg\right) \ar[r] & \IM \left(i_{\MFf_\chi}\right)_{\chi_0} / \left(i_{\MFg}\right)_\chi\left(\tilde{\MCC}_\MFg\right) \ar[r] & \IM \left(i_{\MFf_\chi}\right)_{\chi_0} / \IM \left(i_{\MFg}\right)_\chi \ar[r] & 0.}
\end{equation}
From Lemma \ref{cccodglh}, we deduce the existence of $m_\chi\in\MBN$ such that 
\begin{equation}
\label{didiffjf}
\Char _{\MCO_\MFf[[T]]} \left( \IM \left(i_{\MFf_\chi}\right)_{\chi_0} / \IM \left(i_{\MFg}\right)_\chi \right) = \left(\MFu^{m_\chi}\right),
\end{equation}
with $m_\chi=0$ in case $\mathrm{(i)}$.
Since $\IM \left(i_{\MFg}\right)_\chi / \left(i_{\MFg}\right)_\chi\left(\tilde{\MCC}_\MFg\right) \simeq \left(\MCO_\MFf \widehat{\otimes}_{\MBZ_p} \left( \MCU_{\MFg,\infty} / \MCC_{\MFg,\infty} \right)\right)_\chi$, from (\ref{ddfff}) and (\ref{didiffjf}), we deduce that
\begin{equation}
\label{didiffddjf}
\Char _{\MCO_\MFf[[T]]} \left(\MCO_\MFf \widehat{\otimes}_{\MBZ_p} \left( \MCU_{\MFg,\infty} / \MCC_{\MFg,\infty} \right)\right)_\chi = \MFu^{-m_\chi} \Char _{\MCO_\MFf[[T]]} \left( \IM \left(i_{\MFf_\chi}\right)_{\chi_0} / \left(i_{\MFg}\right)_\chi\left(\tilde{\MCC}_\MFg\right) \right).
\end{equation}
We set $\tilde{\Psi} := \left( \MCO_\MFf \widehat{\otimes}_{\MBZ_p} \overline{\Psi}' \left( K_{\MFf_\chi,\infty} , \MFf_\chi\MFp^\infty \right) \right)_{\chi_0}$.
From (\ref{didiffddjf}) and Lemma \ref{kdugng}, we deduce 
\begin{equation}
\label{didifiiddjf}
\Char _{\MCO_\MFf[[T]]} \left(\MCO_\MFf \widehat{\otimes}_{\MBZ_p} \left( \MCU_{\MFg,\infty} / \MCC_{\MFg,\infty} \right)\right)_\chi = \MFu^{-m_\chi} \Char _{\MCO_\MFf[[T]]} \left( \IM \left(i_{\MFf_\chi}\right)_{\chi_0} / \left(i_{\MFf_\chi}\right)_{\chi_0}\left(\tilde{\Psi}\right) \right).
\end{equation}
Since $\IM \left(i_{\MFf_\chi}\right)_{\chi_0} / \left(i_{\MFf_\chi}\right)_{\chi_0}\left(\tilde{\Psi}\right) \simeq \left( \IM \left(i_{\MFf_\chi}\right) / \left(i_{\MFf_\chi}\right)\left(\tilde{\Psi}\right) \right)_{\chi_0}$
and since $i_{\MFf_\chi}$ is a pseudo-isomorphism, we deduce from (\ref{didifiiddjf}) and Lemma \ref{NIUbhgREeeE} that
\begin{eqnarray}
\label{duvmppsf}
\Char _{\MCO_\MFf[[T]]} \left(\MCO_\MFf \widehat{\otimes}_{\MBZ_p} \left( \MCU_{\MFg,\infty} / \MCC_{\MFg,\infty} \right)\right)_\chi & = & \MFu^{-m_\chi} \Char _{\MCO_\MFf[[T]]} \left( \IM\left(i_{\MFf_\chi}\right) / \left(i_{\MFf_\chi}\right)\left(\tilde{\Psi}\right) \right)_{\chi_0} \nonumber\\
 & = & \MFu^{-m_\chi} \Char _{\MCO_\MFf[[T]]} \left(\MCM\left(G_{\MFf_\chi,\infty},\MCO_\MFf\right) / \MCJ_{\MFf_\chi}\mu\left(\MFf_\chi\right) \right)_{\chi_0} \nonumber\\
 & = & \MFu^{-m_\chi} \Char _{\MCO_\MFf[[T]]} \left(\MCM\left(\GGG,\MCO_\MFf\right) / \chi'_0\left(\MCJ_{\MFf_\chi} \mu\left(\MFf_\chi\right)\right) \right). \nonumber\\
\end{eqnarray}
First we assume that $\MFf_\chi\neq(1)$.
Then $\chi'_0\left(\mu\left(\MFf_\chi\right)\right) \MCM\left(\GGG,\MCO_\MFf\right) \, / \, \chi'_0\left(\MCJ_{\MFf_\chi} \mu\left(\MFf_\chi\right)\right)$ is isomorphic to $\left(\MCO_\MFf \widehat{\otimes}_{\MBZ_p} \mu_{p^\infty} \left(K_{\MFf_\chi,\infty}\right)\right)_{\chi_0}$, hence pseudo-nul since $\mu_{p^\infty} \left(K_{\MFf_\chi,\infty}\right)$ is finite.
Then from (\ref{duvmppsf}) we deduce 
\begin{eqnarray*}
\Char _{\MCO_\MFf[[T]]} \left(\MCO_\MFf \widehat{\otimes}_{\MBZ_p} \left( \MCU_{\MFg,\infty} / \MCC_{\MFg,\infty} \right)\right)_\chi & = & \MFu^{-m_\chi} \Char _{\MCO_\MFf[[T]]} \left(\MCM\left(\GGG,\MCO_\MFf\right) / \chi'_0\left( \mu\left(\MFf_\chi\right) \right) \MCM\left(\GGG,\MCO_\MFf\right) \right)  \\
 & = & \MFu^{-m_\chi} \tilde{\chi}_0\left(\underline{\mu\left(\MFf_\chi\right)}\right)\MCO_\MFf[[T]],
\end{eqnarray*}
and Theorem \ref{theodufnfkffl} follows in this case.
Now assume $\MFf_\chi=(1)$.
Then we expand $\chi'_0$ to the total faction ring of $\MCM\left(G_{(1),\infty} , \MCO_\MFf\right)$ and to the fraction field of $\MCM\left(\GGG , \MCO_\MFf\right)$.
There is $\GGs\in G_\MFg$ such that $\chi(\GGs)\neq1$.
Then $\chi'_0 \left( \mu\left(1\right) \right) \MCM\left( \GGG , \MCO_\MFf \right) \, / \, \chi'_0\left( \MCJ_{(1)}\mu\left(1\right) \right)$ is pseudo-nul, annihilated by $1-\chi(\GGs)$ and $T$.
Since $\chi'_0\left( \MCJ_{(1)}\mu\left(1\right) \right) \subseteq \MCM\left( \GGG , \MCO_\MFf \right)$, we deduce the inclusion $\chi'_0 \left( \mu\left(1\right) \right) \MCM\left( \GGG , \MCO_\MFf \right) \subseteq \MCM\left( \GGG , \MCO_\MFf \right)$ 
and from (\ref{duvmppsf}) we obtain
\begin{equation*}
\Char _{\MCO_\MFf[[T]]} \left(\MCO_\MFf \widehat{\otimes}_{\MBZ_p} \left( \MCU_{\MFg,\infty} / \MCC_{\MFg,\infty} \right)\right)_\chi = \MFu^{-m_\chi} \Char _{\MCO_\MFf[[T]]} \left(\MCM\left(\GGG,\MCO_\MFf\right) / \chi'_0\left(\mu\left(1\right)\right) \MCM\left(\GGG,\MCO_\MFf\right) \right).  
\end{equation*}
$\mathrm{(i)}$ and $\mathrm{(ii)}$ follow immediately in this case. 
\hfill $\square$

\begin{theo}\label{theodufnfkffldeux}
Let $\MFg\in\MCI$ be such that $\MFg\vert\MFf$.
Let $\chi$ be the trivial character on $G_\MFg$.

$\mathrm{(i)}$ If $p\neq2$ or if $w_\MFg = \left\vert \mu(k) \right\vert$, then $\Char _{\MCO_\MFf[[T]]} \left( \MCO_\MFf \widehat{\otimes}_{\MBZ_p} \left(\MCU_{\MFg,\infty} / \MCC_{\MFg,\infty}\right) \right)_\chi$ is generated by $\tilde{\chi}_0 \left( T\underline{\mu(1)} \right)$.

$\mathrm{(ii)}$ If $p=2$, then $\Char _{\MCO_\MFf[[T]]} \left( \MCO_\MFf \widehat{\otimes}_{\MBZ_p} \left(\MCU_{\MFg,\infty} / \MCC_{\MFg,\infty}\right) \right)_\chi$ is generated by $\MFu^{-m_\chi}\tilde{\chi}_0 \left( T\underline{\mu(1)} \right)$, for some $m_\chi\in\MBN$.
\end{theo}

\noindent\textsl{Proof.}
As in the proof of Theorem \ref{theodufnfkffl}, we have
\begin{equation}
\label{duvmppsfdeux}
\Char _{\MCO_\MFf[[T]]} \left(\MCO_\MFf \widehat{\otimes}_{\MBZ_p} \left( \MCU_{\MFg,\infty} / \MCC_{\MFg,\infty} \right)\right)_\chi = \MFu^{-m_\chi} \Char _{\MCO_\MFf[[T]]} \left(\MCM\left(\GGG,\MCO_\MFf\right) / \chi'_0\left(\MCJ_{(1)}\mu\left(1\right)\right) \right),  
\end{equation}
where $m_\chi\in\MBN$ is zero in case $\mathrm{(i)}$.
But $\chi'_0\left(\MCJ_{(1)}\mu\left(1\right)\right) = \chi'_0\left(T\mu\left(1\right)\right) \MCM\left(\GGG,\MCO_\MFf\right)$, and the theorem follows.
\hfill $\square$

\section{Finiteness of invariants and coinvariants.}

For any $\MFh\in\MCI$, we write $\mathrm{L}_{p,\MFh}$ for the $p$-adic $\mathrm{L}$-function of $k$ with modulus $\MFh$, as defined in \cite[II.4.16]{deshalit87}.
It is the map defined on the set of all continuous group morphisms $\xi$ from $\Gal\left(K_{\MFh,\infty}/k\right)$ to $\MBC_p^\times$ (with $\xi\neq1$ if $\MFh=(1)$), such that 
\begin{equation}
\label{defLpxi}
\mathrm{L}_{p,\MFh} \left(\xi\right) = \int \xi(\GGs)^{-1}.\dd\mu(\MFh)(\GGs).
\end{equation}
Let $n\in\MBN$, and let $\chi$ be an irreducible $\MBC_p$-character on $\Gal\left(k\left(\MFh\MFp^n\right)/k\right)$ (with $\chi\neq1$ if $\MFh=(1)$).
We write $F_\chi$ for the subfield of $k\left(\MFh\MFp^n\right)$ fixed by $\Ker(\chi)$, and we write $\chi_{\PR}$ for the character on $\Gal\left(F_\chi/k\right)$ defined by $\chi$.
By inflation we can consider $\chi$ as a group morphism $\Gal\left(K_{\MFh,\infty}/k\right)\rightarrow\MBC_p^\times$, so that the notation $\mathrm{L}_{p,\MFh} \left(\chi\right)$ makes sense.
As in \cite[II.5.2]{deshalit87}, if $n>0$ we set
\begin{equation}
\label{KYnbhj}
\mathrm{L}_{p,\MFh\MFp^n} \left(\chi\right) := \lA \begin{array}{lll} 
\left(1-\chi_{\PR}\left(\MFp,F_\chi/k\right)\right) \mathrm{L}_{p,\MFh} \left(\chi\right) & \text{if} & \text{$\MFp$ is unramified in $F_\chi$,} \\ 
\mathrm{L}_{p,\MFh} \left(\chi\right) & \text{if} & \text{$\MFp$ is ramified in $F_\chi$.} 
\end{array} \right.
\end{equation}

\begin{lm}\label{Lpchizero}
Let $\MFg\notin\{(0),(1)\}$ be an ideal of $\MCO_k$, and let $\chi$ be an irreducible $\MBC_p$-character on $\Gal\left(k(\MFg)/k\right)$.
If $\chi\neq1$ and if none of the prime ideals dividing $\MFg$ are totally split in $F_\chi/k$, then $\mathrm{L}_{p,\MFg} \left(\chi\right)\neq0$.
If $\chi=1$, if $\MFg$ is a power of a prime ideal, and if $\MFp\nmid\MFg$, then $\mathrm{L}_{p,\MFg} \left(\chi\right)\neq0$.
\end{lm}

\noindent\textsl{Proof.}
Let $\MFq\in\{\MFp,\bar{\MFp}\}$ be such that $\GVf_\MFg(1)$ is a unit above $\MFq$.
Let $U\subset k(\MFg)^\times$ be the subgroup generated by the units above $\MFq$.
Let us fix an embedding $\GGi_p:k^{\alg}\hookrightarrow\MBC_p$.
We define the morphism of $k^{\alg}\left[\Gal\left(k(\MFg)/k\right)\right]$-modules below,
\[\ell_p:k^{\alg}\otimes_{\MBZ}U \rightarrow \MBC_p\left[\Gal\left(k(\MFg)/k\right)\right], \quad a\otimes x \mapsto \GGi_p(a) \sum_{\GGs\in\Gal(k(\MFg)/k)} \log_p \left( \GGi_p \left(x^\GGs\right) \right) \GGs^{-1},\]
where $\log_p$ is the $p$-adic logarithm, as defined in \cite[§4]{iwasawa72}.
It is well known that $\Ker\left(\log_p\right)$ is generated by the roots of powers of $p$, so that $\ell_p$ is injective.
Assume $\mathrm{L}_{p,\MFg} \left(\chi\right)=0$.
From \cite[II.5.2, Theorem]{deshalit87}, we deduce that $e_{\chi^{-1}} \ell_p\left( 1\otimes\GVf_\MFg(1) \right)=0$ in $\MBC_p\left[\Gal\left(k(\MFg)/k\right)\right]$.
Then 
\begin{equation}
\label{dhdudi}
e_{\chi^{-1}}\left(1\otimes\GVf_\MFg(1)\right)=0 \quad \text{in} \quad k^{\alg}\otimes_{\MBZ}U,
\end{equation}
where $\chi$ is identified to a group morphism $\Gal\left(k(\MFg)/k\right) \rightarrow k^{\alg}$ via $\GGi_p$.
If $\chi\neq1$, then from \cite[Th\'eor\`eme 10]{robert73} we deduce the existence of a maximal ideal $\MFr$ of $\MCO_k$, unramified in $F_\chi/k$, such that $\MFr\vert\MFg$, and such that $\chi_{\PR}\left(\MFr,F_\chi/k\right)=1$ (hence totally split in $F_\chi/k$).
If $\chi=1$, from (\ref{dhdudi}) we deduce $N_{k(\MFg)/k}\left(\GVf_\MFg(1)\right)\in\mu(k)$.
Then $\MFg$ must be divisible by at least two distinct prime ideals.
\hfill $\square$

\begin{theo}\label{MCUMCCinvcoinvfinite}
For all $n\in\MBN$, the module of $\GGG_n$-invariants and the  module of $\GGG_n$-coinvariants of $\MCU_\infty/\MCC_\infty$ are finite.
\end{theo}

\noindent\textsl{Proof.}
For $n$ large enough, $K_{\MFf,n} / K_n$ is tamely ramified if $p\neq2$, and if $p=2$ the ramification index is $1$ or $2$.
Hence we deduce that the cokernel of the norm maps $\MCU_{\MFf,\infty}\rightarrow\MCU_\infty$ and $\MCU_{\MFf,\infty} / \MCC_{\MFf,\infty} \rightarrow \MCU_\infty / \MCC_\infty$ are annihilated by $2$.
We deduce that
\begin{equation}
\label{OuEttf}
\Char _{\MBZ_p[[T]]} \left( \MCU_\infty / \MCC_\infty \right) \quad\text{divides}\quad 2^a \Char _{\MBZ_p[[T]]} \left( \MCU_{\MFf,\infty} / \MCC_{\MFf,\infty} \right), \quad \text{for some $a\in\MBN$.}
\end{equation}
By (\ref{OuEttf}), we are reduced to prove Theorem \ref{MCUMCCinvcoinvfinite} in the case $K_\infty = K_{\MFf,\infty}$.
By (\ref{chipart}), we only have to show that $\Char_{\MCO_\MFf[[T]]} \left(\MCO_\MFf\widehat{\otimes}_{\MBZ_p}\MCU_{\MFf,\infty} / \MCO_\MFf\widehat{\otimes}_{\MBZ_p}\MCC_{\MFf,\infty}\right)_\chi$ is prime to $\left((1+T)^{p^n}-1\right)$, 
for all $n\in\MBN$, and all irreducible $\MBC_p$-character $\chi$ on $G_\MFf$.
Let $\chi$ be such a character, and let $\GGz\in\mu_{p^\infty}\left(\MBC_p\right)$.
We choose a maximal ideal $\ell$ of $\MCO_k$, prime to $\MFf\MFp$, such that $\chi_{\PR}\left(\ell , F_\chi / k \right)\neq1$ if $\chi\neq1$, and such that $\ell$ is not totally split in $k_1$ (the subfield of $k_\infty$ fixed by $\GGG^p$) if $\chi=1$.
By Theorem \ref{theodufnfkffl} and Theorem \ref{theodufnfkffldeux}, it suffices to prove $\tilde{\chi}_0 \left( \left(1-\GGs_\ell^{-1}\right) \underline{ \mu\left( \MFf_\chi \right) } \right)\vert_{T=\GGz-1}  \neq0$, where $\GGs_\ell := \left(\ell , K_{\MFf,\infty} / k \right)$.
By (\ref{xomftnf}) and by (\ref{socnfyhhfdeux}), we have
\begin{eqnarray}
\label{sdudfjdfldd}
\tilde{\chi}_0 \left( \left(1-\GGs_\ell^{-1}\right) \underline{ \mu\left( \MFf_\chi \right) } \right) \vert_{T=\GGz-1} 
& = & \tilde{\chi}_0 \left( \tilde{\pi}_{\MFf_\chi\ell,\MFf_\chi} \left( \underline{ \mu\left( \MFf_\chi\ell \right) } \right) \right) \vert_{T=\GGz-1} \nonumber \\
& = & \int_\GGG \GGz^{\GGk(\GGs)} . \dd \chi'_{\MFf_\chi\ell} \left( \mu\left( \MFf_\chi\ell \right) \right) (\GGs) ,  
\end{eqnarray}
where $\chi_{\MFf_\chi\ell}$ is the character on $G_{\MFf_\chi\ell}$ defined by $\chi_0$, and where $\GGk:\GGG \rightarrow \MBZ_p$ is the unique morphism of topological groups such that $\GGk(\GGg)=1$.
From (\ref{sdudfjdfldd}) and (\ref{GGsheufg}) we deduce 
\begin{eqnarray}
\label{sdudfjdflDD}
\tilde{\chi}_0 \left( \left(1-\GGs_\ell^{-1}\right) \underline{ \mu\left( \MFf_\chi \right) } \right) \vert_{T=\GGz-1} 
& = & \sum_{g\in G_{\MFf_\chi\ell}} \chi_{\MFf_\chi\ell}(g) \int_\GGG \GGz^{\GGk(\GGs)} . \dd \left(g^{-1}\right)_\ast\mu\left( \MFf_\chi\ell \right) (\GGs).  \nonumber \\
& = & \sum_{g\in G_{\MFf_\chi\ell}} \chi_{\MFf_\chi\ell}(g) \int_{g\GGG} \GGz^{\GGk\left(g^{-1}\GGs\right)} . \dd \mu\left( \MFf_\chi\ell \right) (\GGs) \nonumber \\
& = & \int_{G_{\MFf_\chi\ell,\infty}} \GGz^{\GGk\left(g_\GGs^{-1}\GGs\right)} \chi_{\MFf_\chi\ell}\left(g_\GGs\right) . \dd \mu\left( \MFf_\chi\ell \right) (\GGs),
\end{eqnarray}
where for any $\GGs\in G_{\MFf_\chi\ell,\infty}$, $g_\GGs$ is the image of $\GGs$ through the projection $G_{\MFf_\chi\ell,\infty} \rightarrow G_{\MFf_\chi\ell}$.
We define $\xi : G_{\MFf_\chi\ell,\infty} \rightarrow \MBC_p^\times$, $\GGs \mapsto \GGz^{\GGk\left(g_\GGs^{-1}\GGs\right)} \chi_{\MFf_\chi\ell}\left(g_\GGs\right)$.
Then $\xi$ is a group morphism, and if $n\in\MBN$ is such that $\GGz^{p^n}=1$, then $\xi$ defines an irreducible $\MBC_p$-character on $G_{\MFf_\chi\ell,n} := \Gal\left(K_{\MFf_\chi\ell,n}/k\right)$.
Let $\MFg$ be the conductor of $F_\xi$.
Since the restriction of $\xi$ to $G_{\MFf_\chi\ell}\hookrightarrow G_{\MFf_\chi\ell,n}$ is $\chi_{\MFf_\chi\ell}$, we deduce that there is $m\in\MBN$ such that $\MFg = \MFf_\chi \MFp^m$, and from (\ref{KYnbhj}) we deduce that 
\begin{equation}
\label{kHiyHJhGgjhJY}
\mathrm{L}_{p,\MFg\ell} \left( \xi^{-1} \right) \, = \, \mathrm{L}_{p,\MFf_\chi\ell} \left( \xi^{-1} \right).
\end{equation}
Then from (\ref{sdudfjdflDD}) and (\ref{defLpxi}) we deduce
\begin{equation}
\label{wpyafjf}
\left(1-\tilde{\chi}_0\left(\GGs_\ell^{-1}\right)\right) \tilde{\chi}_0 \left( \underline{ \mu\left( \MFf_\chi \right) } \right) \vert_{T=\GGz-1} = \int_{G_{\MFf_\chi\ell,\infty}} \xi(\GGs) . \dd \mu\left( \MFf_\chi\ell \right) (\GGs) = \mathrm{L}_{p,\MFf_\chi\ell} \left(\xi^{-1}\right).
\end{equation}
If $\chi\neq1$, then $\chi_{\PR}\left(\ell,F_\chi/k\right)\neq1$ implies that $\ell$ is not totally split in $F_\xi/k$.
If $\chi=1$ and $\GGz\neq1$, then $k_1\subseteq F_\xi$ and $\ell$ is not totally split in $F_\xi/k$.
If $\chi=1$ and $\GGz=1$, then $\xi=1$ and $\MFg=(1)$.
From (\ref{wpyafjf}), (\ref{kHiyHJhGgjhJY}), and Lemma \ref{Lpchizero}, we deduce 
\begin{equation*}
\left(1-\tilde{\chi}_0\left(\GGs_\ell^{-1}\right)\right) \tilde{\chi}_0 \left( \underline{ \mu\left( \MFf_\chi \right) } \right) \vert_{T=\GGz-1} = \mathrm{L}_{p,\MFg\ell} \left(\xi^{-1}\right) \neq 0.
\end{equation*}
\hfill $\square$


\bibliographystyle{amsplain}

\end{document}